\documentclass{ajour}
\usepackage{epsfig}

\usepackage{graphicx}
\usepackage{amsmath}
\usepackage{amsfonts}
\usepackage{latexsym}
\begin{document}
\authorrunninghead{Giuseppe Melfi}
\titlerunninghead{On simulataneous binary expansions of $n$ and $n^2$}
\title{On simultaneous binary expansions of $n$ and $n^2$}

\author{Giuseppe Melfi}

\affil{Universit\'e de Neuch\^atel\\
Groupe de Statistique \\
Espace de l'Europe 4, CH--2002 Neuch\^atel, Switzerland}

\email{Giuseppe.Melfi@unine.ch
}

\abstract{A new family of sequences is proposed.
An example of sequence of this family is more accurately studied.
This sequence is composed by the integers $n$ for which the sum of binary 
digits is equal to the sum of binary digits of $n^2$. 
Some structure and asymptotic 
properties are proved and a conjecture about its counting function is
discussed.}


\begin{article}


\keywords{Sequences and sets. Binary expansion.
Asymptotic behaviour.}

\section{Introduction}
This paper arises from some natural questions about digits of numbers. 
We are interested in numbers $n$
such that $n$ and $n^m$ have a certain relation involving
sums of digits.
Let $B(n)$ be the sum of digits of the positive integer $n$
written on base $2$.
This function represents the numbers of ones in the binary
expansion of $n$, or from the code theory point of view, 
the number of nonzero digits in the bits string representing $n$,
i.e., the so-called `Hamming weight' of $n$. 
It is obvious that $B(n)=B(2n)$, but
for a prime $p>2$  a  relation between $B(n)$
and $B(pn)$ is less trivial
\cite{coquet1,delange,dumont,grabner,schmid}.  

In this paper we are interested in comparing $B(n)$ with $B(n^m)$.
Stolarsky \cite{stolarsky} proved some inequalities for the
functions $r_m(n)=B(n^m)/B(n)$.
Lindstr\"om \cite{lindstrom} proved that $\limsup_{n\rightarrow\infty}
{B(n^m)}/\log_2n=m$,
for $m\ge2$. 

One can naturally define a family of positive integer sequences. 

\begin{definition}
Let $k\ge2$, $l\ge1$, $m\ge2$ be positive integers.
We say that a positive integer $n$ satisfies the property $(k,l,m)$ if the 
sum of its digits
in its expansion in base $k$ is $l$ times the sum of the digits of the
expansion in base $k$ of $n^m$.
\end{definition}  

For every triplet $(k,l,m)$ we have a sequence made up of the positive
integers of the type $(k,l,m)$. 

The simplest case is $(k,l,m)=(2,1,2)$, which corresponds to the
positive integers $n$ for which the numbers of ones in their binary
expansion is equal to the number of ones in $n^2$.

This paper is consacrated to the study of the $(2,1,2)$-numbers.
The list shows several interesting facts (see Table~1). 
The distribution is not
regular.  A huge amount of questions, most of which of
elementary nature, can be raised. 

In despite of its elementary definition, this sequence surprisingly
does not appear in literature. Recently, proposed by the author,
it appeared on \cite{Sloaneonline}.

Several questions, concerning both the structure 
properties and asymptotic behaviour, can be raised. 
Is there a necessary and sufficient 
condition to assure that a number is of type (2,1,2)? What is the asymptotic
behaviour of the counting function of (2,1,2)-numbers? 

The irregularity of distribution
does not suggest a clear answer to these questions.

In Section~\ref{arit} we provide some sufficient conditions in order
that a number is of type (2,1,2). In particular we explicitely 
provide several infinite sets of (2,1,2)-numbers, but all these sets
are quite thin, and the problem of an exhaustive answer does not appear
easy at first sight.
.

Let $p(n)$ be the number of (2,1,2)-numbers which does not exceed $n$.

\begin{conjecture}\label{conjecture}
Let $\displaystyle p(n)=\sum_{m<n\atop \tiny m \mbox{of type (2,1,2)}}1\mbox{ }$ 
be the counting function of iso\-square numbers. There exists a continuous
function $F(x)$, periodic of period $1$, such that 
$$p(n)= \frac{n^\alpha}{\log n}F\left(\frac {\log n}{\log2}\right)+R(n),$$
where $\alpha=\log{1.6875}/\log2\simeq0.75488750$, and 
$R(n)=o(n^\alpha/\log n)$.

Further $F(n)$ is nowhere differentiable.
\end{conjecture}

Sequences constructed in Section~\ref{arit} yield $p(n)\gg\log n$.
In Section~\ref{asym} we prove a lower bound
$p(n)\gg n^{0.025}$. Conjecture~\ref{conjecture} is suggested by a
more detailed observation of the list, as well as some similarities
with certain functions studied by Boyd {\em et al.} \cite{boyd}.

A discussion on Conjecture~\ref{conjecture} 
will be developed in Section~\ref{conj}.
Our arguments appear to  be confirmed by experimental results 
done for $n<10^8$.

\section{Notations}
If the binary expansion of $n$ is $c_1c_2\dots c_k$, ($c_i\in\{0,1\}$) we will
write $n=(c_1c_2\dots c_k)$. The first digit may be $0$ so, for example,
we allow the notation $3=(011)$.
We will also denote $\overbrace{1\dots1}^{\tiny
k\mbox { times}}$ as $(1_{(k)})$, so for example,
$$
(\overbrace{11\dots11}^{\tiny k\mbox { times}}
\overbrace{00\dots00}^{\tiny h\mbox { times}}
\overbrace{11\dots11}^{\tiny l\mbox  { times}}010
)=(1_{(k)}0_{(h)}1_{(l)}010).
$$
This allows to perform arithmetical operations in a more compact manner. 
For example, if
$k>h+h'$, one has $(1_{(k)})-(1_{(h)}0_{(h')})=(1_{(k-h-h')}0_{(h)}1_{(h')})$.

Let $B(n)$ be the number of 1's in the binary expansion of $n$. 
This function, known also as the Hamming weight of $n$, 
is used mainly in the theory of algorithms, namely for the study 
of computational aspects and complexity. 
So, $n$ is of type (2,1,2) if $B(n)=B(n^2)$. 

Let $c\in\{0,1\}$ a binary digit. 
We will use the notation $c'$ to indicate the other
digit: $c'=1-c$. This notation will be very useful in computing the Hamming
weight of a number. We will use the property that
if $n=(c_1c_2\dots c_k)$, then $2^k-n-1=(c_1'c_2'\dots c_k')$, and
$B(n)=k-B(2^k-n-1)$.

\section{Arithmetics and structure properties}\label{arit}

There are infinitely many (2,1,2)-numbers. 
Note that if $n$ is an (2,1,2)-number, then $2n$ is an (2,1,2)-number,
and that if $n$ is an even (2,1,2)-number, then $n/2$ is a
(2,1,2)-number. 
There is no direct dependence between $2n+1$  or $2n-1$  and $n$, 
classical arguments for the study of asymptotic properties of the counting function
(see \cite{tenenbaum}) cannot be applied.
In the following remarks we show some sets of (2,1,2)-numbers.

\begin{remark}\label{base}
For every $k>1$, the number $n_k=2^k-1$ is of type (2,1,2).
\end{remark}
In fact, $B(n_k)=k$, and
$$(2^k-1)^2=2^{2k}-2^{k+1}+1=(\overbrace
{11\dots11}^{\tiny k-1\mbox{ times}}
\overbrace
{00\dots00}^{\tiny k\mbox{ times}}1),$$
so $n_k^2$ also contains $k$ times the digit $1$.

Table 1 contains 4-tuples of consecutive integers of (2,1,2)-numbers.
After the 4-tuple $(1,2,3,4)$ the second one is $(316,317,318,319)$. 

\begin{remark}
There are infinitely many 4-tuples of consecutive integers composed by
(2,1,2)-numbers.
\end{remark}

In fact it is an exercice to show that for every $k\ge9$ and 
$n=2^k-2^{k-2}-2^{k-3}-4$, the numbers $n$, $n+1$, $n+2$ and
$n+3$ are all of type (2,1,2).  

The following proposition shows that the set of (2,1,2)-numbers presents
a certain gap structure.

\begin{proposition}
There are infinitely many $n$ such that the interval
$]n,n+n^\frac12[$ does not contain
any (2,1,2)-number.
\end{proposition}

\begin{proof}
Let $n=2^{2k}=(10_{(2k)})$. Each $m\in ]n,n+n^\frac12[$ is of the form
$n+r$ with $r<n^\frac12$. In its binary espansion $m$
is of the kind
$(10_{(k)}c_1c_2\dots c_k)$. Here $c_i\in\{0,1\}$
are binary digits and $B((c_1c_2\dots c_k))=B(r)\ge1$. Let
$r^2=(q_1q_2\dots q_{2k})$.  
We have again $B((q_1q_2\dots q_{2k}))\ge1$. Hence
$$\begin{array}{rl}
m^2&=\{2^{2k}+( c_1c_2\dots c_k)\}^2\\
&\\
   &=2^{4k}+2^{2k+1}( c_1c_2\dots c_k)+( c_1c_2\dots c_k)^2 \\
&\\
   &=(10_{(k-1)}c_1c_2\dots c_k0q_1q_2\dots q_{2k})
\end{array}
$$
Therefore $B(m^2)=1+B(r)+B(r^2)>1+B(r)=B(m)$.
\end{proof}

The preceding construction cannot be improved. One can easily prove that
if $n>5$ is odd, then there exists an (2,1,2)-number between $2^n$ and
$2^n+4\cdot2^\frac n2$. Namely if $n=2m+1$ such a number is 
$2^{2m+1}+2^{m+2}-1$. The proof by check digits in column operations is
straightforward.  

\section{A lower bound for the counting function}\label{asym}

We begin this section with some preliminary lemmata.

\begin{lemma}\label{lem1}
If $n<2^\nu$, then $B(n(2^\nu-1))=\nu$.
\end{lemma}

\begin{proof}
We assume that $n$ is odd.
Let $n=(c_1c_2\dots c_k)$, with $c_1=c_k=1$. 
Since $n<2^\nu$ we have that $k\le \nu$.
We have that $n(2^\nu-1)=(c_1c_2\dots c_k0_{(\nu)})-(c_1c_2\dots c_k)$. 

Hence
$n(2^\nu-1)=(c_1c_2\dots c_{k-1}c_k'1_{(\nu-k)}c_1'c_2'\dots c_{k-1}'c_k)$.
Therefore
$$
\begin{array}{rl}
B(n(2^\nu-1))&=B((c_1c_2\dots c_{k-1}c_k'1_{(\nu-k)}c_1'c_2'\dots c_{k-1}'c_k))\\
&\\
&=B((c_1c_2\dots c_{k-1}c_k1_{(\nu-k)}c_1'c_2'\dots c_{k-1}'c_k'))\\
&\\
&=B(n)+(\nu-k)+(k-B(n))\\
&\\
&=\nu.
\end{array}
$$
If $n$ is even, $n'=n/2^h$ is an odd integer for a certain $h$, and 
$n'<2^\nu$. Hence $B(n'(2^\nu-1))=\nu$ and 
$B(2^hn'(2^\nu-1))=B(n'(2^\nu-1))=\nu$, so $B(n(2^\nu-1))=\nu$.
\end{proof}

\begin{lemma}\label{lem3}
Let $n\in\mathbb N$ and let $\nu$ be such that $n<2^{\nu-1}$. Then
$$ B(2^\nu n+1)=B(n)+1 \ \ \mbox{ and } B( (2^\nu n+1)^2)=B(n^2)+B(n)+1.$$
\end{lemma}
\begin{proof} 
The proof is straightforward.\end{proof}

\begin{lemma}\label{lem4}
Let $n=(c_1c_2\dots c_k)$, $m=(d_1d_2\dots d_h)$ 
odd positive integers. If 
$\nu\ge\max\{2h-1,h+k+1\}$, we have
$$B(n2^\nu-m)=B(n)-B(m)+\nu$$
and
$$B((n2^\nu-m)^2)=B(n^2)+B(m^2)-B(mn)+\nu-1.$$
\end{lemma} 

\begin{proof}
Let $n=(c_1c_2\dots c_k)$ and $ m=(d_1d_2\dots d_h)$. If $h\le\nu$ we have
$n2^\nu-m=(c_1c_2\dots c_{k-1}c_k'1_{(\nu-h)}d_1'd_2'\dots d_{h-1}'d_h)$.
Since $m$ and $n$ are odd, $c_k'+d_h=c_k+d_h'$, so
$B(n2^\nu-m)=B((c_1c_2\dots c_k1_{(\nu-h)}d_1'd_2'\dots d_h'))=B(n)+\nu
-B(m).$

Let $n^2=(q_1q_2\dots q_{2k})$ and $m^2=(r_1r_2\dots r_{2h})$. Let
$mn=(s_1s_2\dots s_{k+h})$. We have 
$(n2^\nu-m)^2=(n^22^{2\nu}+m^2)-nm2^{\nu+1}$, and if $2h\le\nu+1$
and $k+h\le\nu-1$, we have
$$(n2^\nu-m)^2=(q_1\dots q_{2k-1}q_{2k}'1_{(\nu-k-h-1)}s_1'\dots
s_{k+h-1}'s_{k+h}0_{(\nu+1-2h)}r_1\dots r_{2h}),$$
hence $B((n2^\nu-m)^2)=B(n^2)+B(m^2)-B(mn)+\nu-1$.
\end{proof}

\begin{corollary}
\label{coroold}
If $B(n^2)=2B(n)-1$, and $\nu\ge k+2$,
then $2^\nu n-1$ is of type (2,1,2).
\end{corollary}

\begin{corollary}\label{coro3}
Let $n$ an odd positive integer. Let $m=2^{h}-1$ with $n<m$. If
$\nu\ge 2h+1$ then
$$B(n2^\nu-m)=B(n)+\nu-h$$
and
$$B((n2^\nu-m)^2)=B(n^2)+\nu-1.$$
\end{corollary}
\begin{proof}
These statements are an easy consequence of Lemma~\ref{lem4}, of 
Lemma~\ref{lem1}, and  of Remark~\ref{base}.
\end{proof}

\begin{lemma}\label{lem5}
Let $n=(c_1c_2\dots c_k)$ be an odd positive integer. 
Let $B(n^2)\ge2B(n)+1$.
There exist $\nu$ and $h\in\mathbb N$ such that for
$n'= n2^\nu-(2^{h}-1)$ we have
$$B(n'^2)=2B(n')-1.$$
\end{lemma}

\begin{proof}
Let $h=k+1$. We have $n<2^h-1$.
Let $\nu=B(n^2)-2B(n)+2h$. We have $\nu=2h+a$ with $a\ge1$. 
Remark that $\nu<4k+2$.
The hypotheses of Corollary~\ref{coro3} are satisfied, so
$B(n')=B(n)+\nu-h$ and
$$\begin{array}{rl}
B(n'^2)&=B(n^2)+\nu-1\\
&\\
&=2B(n)+\nu-2h+\nu-1\\
&\\
&=2B(n)+2\nu-2h-1\\
&\\
&=2B(n')-1.
\end{array}
$$
\end{proof}

\begin{lemma}\label{lem6}
Let $n=(c_1c_2\dots c_k)>1$ be a positive integer. Let
$n_0=(c_1c_2\dots c_k0_{(k)}10_{(2k+1)}1)$. Then
$$B(n_0^2)>2B(n_0)+1.$$
\end{lemma}

\begin{proof}
It suffices to apply twice Lemma~\ref{lem3}. Remark that $n_0\ll n^4$.
\end{proof}

\begin{theorem}
Let $p(n)$ be the counting function of (2,1,2)-numbers. We have
$$p(n)\gg n^{0.025}.$$
\end{theorem}

\begin{proof}
Let $n=(c_1c_2\dots c_k)$ be an odd positive integer. 
We will show that for a constant $A$ not depending on $n$,
it is possible to construct a set of $n$ distinct (2,1,2)-numbers
not exceeding $An^{40}$.

Let consider the $n$ integers $n_i=2^k+i$,
for $i=1,\dots n$. We have obviously  $n_i<4n$. 

For every $i$, by Lemma~\ref{lem6} it exists
$n_{0,i}\ll n_i^4$ whose first $k+1$ digits are the same as those of $n_i$
such that $B(n_{0,i}^2)>2B(n_{0,i})+1$.

By Lemma~\ref{lem5} it exists $n_{0,i}'\ll n_{0,i}^5$ such that
the first $k$ binary digits of $n_{0,i}'$ are again those of $n$ and
such that $B( n_{0,i}'^2)=2B(n_{0,i}')-1$.

Finally, by Corollary~\ref{coroold}, it exists $n_{0,i}''\ll n_{0,i}'^2$,
whose the first binary digits are the same as for $n$ and such that
$B(n_{0,i}''^2)=B(n_{0,i}'')$.

We have $n_{0,i}''\ll n_{0,i}'^2$ $\ll (n_{0,i}^5)^2\ll$
$((n_i^4)^5)^2\ll n^{40}$.
\end{proof}

\section{A probabilistic approach}\label{conj}

Let $k$ be a sufficiently large positive integer. Let $n$ be such that 
$2^k\le n<2^{k+1}$. In a binary base, these numbers are made up of a
`1' digit followed by $k$ binary digits $0$ or $1$. So $1\le B(n) \le k+1$.
Let us consider $n^2$. 
We have $2^{2k}\le n^2<2^{2k+2}$, so its binary expansion
contains a `1' digit followed by $2k$ or $2k+1$ binary digits $0$ or $1$.

In this section we estimate the asymptotic behaviour of $p(n)$
under a suitable assumption. 

We consider $B(n)$ and $B(n^2)$ as random variables. 
Clearly $B(n)-1$ follow a binomial random distribution of type
$b(k,1/2)$. Schmid \cite{schmid} studied joint distribution
of $B(p_in)$ for distinct odd integer $p_i$. Here we consider the 
joint distribution of $B(n)$ and $B(n^2)$.
We assume that for sufficiently large $k$, 
$B(n^2)-1$ tends to follow a binomial random distribution of type
$b(2k,1/2)$ if $2^{2k}\le n^2<2^{2k+1}$ and a binomial random 
distribution of type
$b(2k+1,1/2)$ if $2^{2k+1}\le n^2<2^{2k+2}$. We assume that $B(n)$ and 
$B(n^2)$ are independent realizations of these random variables.

It is clear that for very small values of $B(n)$, the numerical value of
$B(n^2)$ is also relatively small, since $B(n^2)\le B(n)^2$, 
so these variables are not completely 
independent. But for $B(n)>\sqrt{\log n}$ this phenomenon tends to disappear, 
and the 
preceding assumption can be taken in account to an asymptotic behaviour
estimate.

Hence 

$$\Pr(n \mbox{ of type }(2,1,2) \mbox{ and } B(n)=l)=\frac{
\left(\begin{array}{c}2k\\l\end{array}\right)+
\left(\begin{array}{c}2k+1\\l\end{array}\right)}
{3\cdot2^{2k}}$$

This suggests that
$p(n)\sim n^\alpha$
with
$$\alpha=-2+\frac1{\log2}\lim_{k\rightarrow\infty}
\frac1k\log\sum_{l=0}^k\left(\begin{array}{c}k\\l\end{array}\right)
\left(\begin{array}{c}2k\\l\end{array}\right)=\frac{\log{27/16}}{\log2}\simeq 
0.75488.$$


The least square method applied to (2,1,2)-numbers not exceeding $10^8$ gives
the value $\alpha\simeq0.73$. The difference is due to the fact that there is
a small effect of correlation between $B(n)$ and $B(n^2)$ for $n$ in a 
neighbourhood of a power of $2$. 
Hence we conjecture that $p(n)\asymp n^\alpha/\log n$.

However, a plot of $p(n)\log n/n^\alpha$ shows more complex details
(see Figure~1). It seems that
the limiting function is not derivable. Effectively, as shown by 
Delange \cite{delange} and by Coquet \cite{coquet1}, and more recently by
Tenenbaum \cite{tenenbaum}, if $f(n)$ is a function
related to the binary expansion of $n$, often one has
that $F(n)=\sum_{k\le n}{f(k)}$ has an expression in which periodic functions,
often nowhere differentiable, are involved.
These properties probably hold in our case, but a direct approach as shown in 
\cite{tenenbaum} does not appear possible.
This remark, joint with the observation of Figure 1, justifies
Conjecture~\ref{conjecture}.

\section*{Acknowledgments}
I am grate to Chu Wen-Chang and to Roberto Avanzi  for 
their useful discussions. I thank Alina Matei for her help in 
the early phases of the computational approach.

This work has been partly realized during my stay in Cetraro (Italy)
with the contribution of C.I.M.E. I am grate to Carlo Viola
for his strong encouragement and support.

\section*{Appendix}

\vspace{-3mm}
\input plot2

\noindent
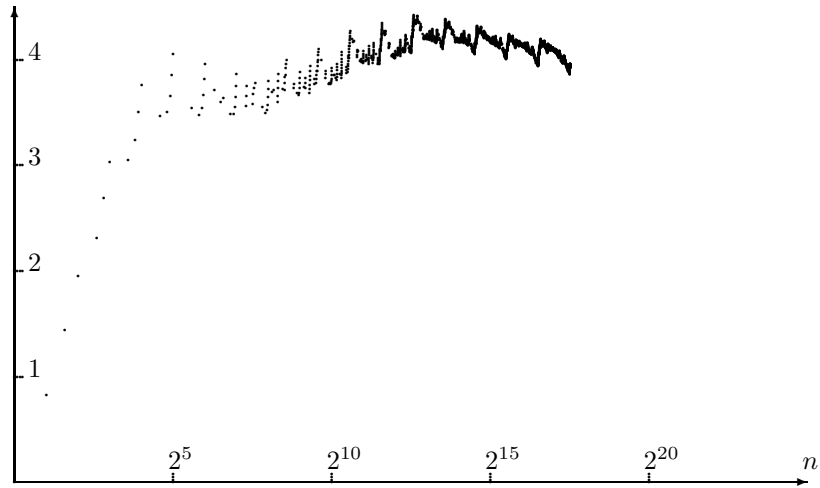
Figure 1. The plot of $p(n)\log n/n^\alpha$. A logarithmic scale is used in abscissae.

{\tiny 
$$
\begin{array}{rrrrrrrrrrrrrrr}
1&248&702&1272&1951&2812&3560&4594&5624&6124&7647&8701&10142&11642&12254\\
2&252&703&1274&1984&2814&3572&4596&5628&6126&7670&8702&10144&11648&12255\\
3&254&728&1276&2016&2815&3578&4600&5630&6127&7671&8703&10176&11676&12268\\
4&255&730&1277&2032&2898&3581&4601&5631&6134&7675&8750&10192&11680&12270\\
6&256&732&1278&2040&2910&3584&4602&5796&6135&7680&8928&10208&11684&12271\\
7&279&735&1279&2044&2912&3630&4604&5799&6139&7800&8958&10216&11694&12278\\
8&287&748&1404&2046&2919&3806&4605&5807&6144&7804&8959&10224&11704&12279\\
12&314&750&1406&2047&2920&3835&4606&5811&6271&7805&9160&10228&11708&12283\\
14&316&751&1407&2048&2921&3840&4607&5820&6395&7806&9168&10232&11712&12288\\
15&317&758&1449&2142&2926&3900&4803&5821&6396&7870&9184&10234&11726&12542\\
16&318&759&1455&2159&2927&3902&4859&5824&6398&7871&9188&10236&11730&12543\\
24&319&763&1456&2173&2928&3903&4860&5838&6399&7903&9192&10237&11760&12639\\
28&351&768&1460&2174&2940&3935&4862&5840&6477&7936&9200&10238&11822&12790\\
30&364&815&1463&2175&2975&3968&4863&5842&6479&8031&9202&10239&11900&12792\\
31&365&890&1464&2232&2987&4032&4911&5847&6518&8047&9204&10462&11901&12795\\
32&366&893&1470&2290&2990&4064&4991&5852&6520&8064&9208&10494&11903&12796\\
48&374&896&1495&2292&2992&4080&5020&5854&6523&8128&9209&10495&11942&12798\\
56&375&960&1496&2296&3000&4088&5024&5856&6557&8160&9210&10740&11948&12799\\
60&379&975&1500&2297&3002&4092&5055&5863&6559&8176&9212&10746&11960&12954\\
62&384&992&1501&2298&3004&4094&5056&5865&6588&8184&9213&10749&11968&12958\\
63&445&1008&1502&2300&3007&4095&5069&5880&6623&8188&9214&10837&11999&13036\\
64&448&1016&1516&2301&3032&4096&5071&5911&6638&8190&9215&11006&12000&13040\\
79&480&1020&1518&2302&3036&4191&5072&5950&6644&8191&9469&11007&12008&13046\\
91&496&1022&1519&2303&3038&4207&5088&5971&6647&8192&9606&11230&12015&13051\\
96&504&1023&1526&2430&3039&4223&5096&5974&6650&8382&9718&11231&12016&13114\\
112&508&1024&1527&2431&3052&4253&5104&5980&6653&8414&9720&11232&12028&13118\\
120&510&1071&1531&2510&3054&4284&5108&5984&6831&8415&9723&11247&12128&13176\\
124&511&1087&1536&2512&3055&4318&5112&6000&6871&8445&9724&11248&12131&13233\\
126&512&1116&1599&2528&3062&4345&5114&6004&6879&8446&9726&11256&12144&13246\\
127&558&1145&1630&2536&3063&4346&5116&6008&7101&8447&9727&11260&12152&13276\\
128&573&1146&1647&2544&3067&4348&5117&6014&7119&8506&9822&11262&12154&13288\\
157&574&1148&1661&2548&3072&4349&5118&6064&7120&8568&9982&11263&12156&13294\\
158&575&1149&1780&2552&3198&4350&5119&6072&7144&8636&9983&11591&12159&13300\\
159&628&1150&1786&2554&3199&4351&5231&6076&7156&8689&10035&11592&12208&13303\\
182&632&1151&1789&2556&3259&4375&5247&6077&7162&8690&10040&11597&12216&13306\\
183&634&1215&1792&2557&3260&4464&5370&6078&7165&8692&10048&11598&12220&13309\\
187&636&1255&1815&2558&3294&4479&5373&6104&7168&8696&10079&11614&12222&13591\\
192&637&1256&1903&2559&3319&4580&5503&6108&7260&8697&10110&11622&12223&13662\\
224&638&1264&1920&2685&3322&4584&5615&6110&7455&8698&10112&11639&12248&13739\\
240&639&1268&1950&2808&3325&4592&5616&6111&7612&8700&10138&11640&12252&13742\\
\end{array}
$$
}

\begin{center}
Table 1. The $(2,1,2)$-numbers not exceeding 13742.
\end{center}

\end{article}
\end{document}